\documentclass[a4paper,11pt,reqno]{amsart}

\usepackage{latexsym,dsfont,amssymb,amsmath,amsthm}

\chardef\bslash=`\\ % p. 424, TeXbook

\newtheorem{theorem}{Theorem}
\newtheorem{corollary}[theorem]{Corollary}
\newtheorem{lemma}[theorem]{Lemma}

\theoremstyle{remark}

\theoremstyle{definition}

\newcommand\bp{\begin{proof}}
\newcommand\ep{\end{proof}}

\newcommand{\Z}{\mathbb Z}

\newcommand\T{\mathbb T}

\newcommand\Ind{\operatorname{Ind}}

\newcommand{\Ad}{\operatorname{Ad}}

\newcommand{\supp}{\operatorname{supp}}

\newcommand\enu[1]{\smallskip\newline\makebox[5mm][l]{\rm(#1)}}

\begin{document}

\title[Traces on crossed products]
{Traces on crossed products}

\author[S. Neshveyev]{Sergey Neshveyev}
\address{Department of Mathematics, University of Oslo,
P.O. Box 1053 Blindern, N-0316 Oslo, Norway}

\email{sergeyn@math.uio.no}

\thanks{Partially supported by the Research Council of Norway}

\begin{abstract}
We give a description of traces on $C(X)\rtimes G$ in terms of measurable fields of traces on the C$^*$-algebras of the stabilizers of the action of $G$ on $X$.
\end{abstract}

\date{October 4, 2010}

\maketitle

The present note is motivated by a study of KMS-states on crossed products $C(X)\rtimes S$, where~$S$ is a group-like object (group, semigroup, Hecke algebra,...) and the dynamics is given by a one-parameter subgroup of the dual action. One of the interesting examples is the system of Bost and Connes~\cite{bos-con}. A~standard strategy, see e.g.~\cite{LLNlat}, is to consider the measure on $X$ defined by a KMS-state and show that the KMS-condition forces the set of points with nontrivial stabilizers to have measure zero. This usually implies that the state factors through the conditional expectation onto $C(X)$, which reduces the study of KMS-states to a dynamical systems problem. But in general it is clear that stabilizers fixed by the dynamics can carry nontrivial traces. Our goal is to describe precisely how this happens in the case of crossed products by group actions. This description is particularly transparent for abelian groups.

\medskip

We start with the following observation, which in a way goes back to \cite[Appendix]{MvN}. Recall that the centralizer $A_\varphi$ of a state $\varphi$ on a C$^*$-algebra~$A$ is the set of elements $a\in A$ such that $\varphi(ab)=\varphi(ba)$ for all $b\in A$.

\begin{lemma} \label{lcent}
For any state $\varphi$ on a unital C$^*$-algebra $A$, there exists a unique state $\Phi$ on $A_\varphi^{op}\otimes_{max}A$ such that $\Phi(a\otimes b)=\varphi(ab)$ for all $a\in A_\varphi$ and $b\in A$.
\end{lemma}

\bp We may assume that $A\subset B(H)$ and $\varphi$ is defined by a cyclic vector~$\xi\in H$. Assume first that~$\xi$ is separating for $A''$. Let $J$ be the corresponding modular conjugation. Define a representation $\pi$ of $A_\varphi^{op}\otimes_{max}A$ on $H$ by
$\pi(a\otimes b)=Ja^*Jb=bJa^*J$. If $a\in A_\varphi$, then~$a$ commutes with the modular operator, hence $Ja\xi=a^*\xi$.
Therefore $\pi(a\otimes b)\xi=ba\xi$, so that $\Phi:=(\pi(\cdot),\xi,\xi)$ is the required state.

\smallskip

In the general case we will show that a representation $\pi$ of $A_\varphi^{op}\otimes_{max}A$ such that $\pi(a\otimes b)\xi=ba\xi$ always exists. For $a\in A_\varphi$ and $b,c\in A$ we have
$$
(ba\xi,c\xi)=(b\xi,ca^*\xi).
$$
It follows that for every $a\in A_\varphi$ there exists a well-defined operator $\rho(a)$ on $A\xi$ such that $\rho(a)b\xi=ba\xi$, and then $\rho$ is a representation of $A^{op}_\varphi$ on $A\xi$. Since $\rho(u)$ is unitary for unitary $u$, this is a representation by bounded operators, so it extends to a representation of $A^{op}_\varphi$ on $H$. Its image commutes with $A$, so we can define a representation $\pi$ of $A_\varphi^{op}\otimes_{max}A$ on $H$ by $\pi(a\otimes b)=\rho(a)b$. Then $\pi(a\otimes b)\xi=ba\xi$.
\ep

Assume now that a countable group $G$ acts by homeomorphisms on a metrizable compact space~$X$. Denote by $u_g\in C(X)\rtimes G$ the canonical unitaries implementing the action: $u_gfu^*_g=f(g^{-1}\cdot)$. Consider also the full group C$^*$-algebra $C^*(G)$ of $G$ with generators $\lambda_g$, $g\in G$. Denote by $j\colon C^*(G)\to C(X)\rtimes G$ the canonical homomorphism mapping $\lambda_g$ into $u_g$. Assume $\varphi$ is a state on $C(X)\rtimes G$ with centralizer containing~$C(X)$. By Lemma~\ref{lcent} we can define a state $\Phi$ on $C(X)\otimes C^*(G)$ by
\begin{equation}\label{etrace}
\Phi(f\otimes a)=\varphi(fj(a))\ \ \hbox{for}\ \ f\in C(X)\ \ \hbox{and}\ \ a\in C^*(G).
\end{equation}
Disintegrating $\Phi$ with respect to $C(X)$ we get a probability measure $\nu$ on $X$ and states $\varphi_x$ on~$C^*(G)$ such that $\Phi=\int^\oplus_X\varphi_xd\nu(x)$. Then
\begin{equation*}
\varphi(fu_g)=\Phi(f\otimes\lambda_g)=\int_Xf(x)\varphi_x(\lambda_g)d\nu(x)\ \ \hbox{for}\ \ f\in C(X)\ \ \hbox{and}\ \ g\in G.
\end{equation*}
Denote by $G_x$ the stabilizer of $x\in X$.

\begin{theorem} \label{tmain}
Identity~\eqref{etrace} defines a bijection between states $\varphi$ on $C(X)\rtimes G$ with centralizer containing~$C(X)$ and states $\Phi=\int^\oplus_X\varphi_xd\nu(x)$ on $C(X)\otimes C^*(G)$ such that $\varphi_x(\lambda_g)=0$ for $\nu$-a.e.~$x\in X$ and every $g\notin G_x$.
\end{theorem}

\bp Assume $\varphi$ is a state with centralizer containing $C(X)$, and $\Phi$ is the state defined by~\eqref{etrace}. Fix an element $g\in G$, $g\ne e$. Assume a point $x_0\in X$ is not fixed by $g$. Choose an open neighbourhood~$U$ of~$x_0$ such that $U\cap gU=\emptyset$. Then for any continuous functions $f_1$ and~$f_2$ with supports in $U$ we have
$
\varphi(f_1f_2u_g)=\varphi(f_2u_gf_1)=\varphi(f_2f_1(g^{-1}\cdot)u_g)=0,
$
that is,
$$
\int_Xf_1(x)f_2(x)\varphi_x(\lambda_g)d\nu(x)=0.
$$
It follows that $\varphi_x(\lambda_g)=0$
for $\nu$-a.e.~$x\in U$. Since this is true for any sufficiently small neigbourhood~$U$ of a point not fixed by~$g$, we conclude that $\varphi_x(\lambda_g)=0$ for $\nu$-a.e.~$x$ not fixed by $g$. Since $G$ is countable, this proves the theorem in one direction.

\smallskip

Conversely, assume $\Phi=\int^\oplus_X\varphi_xd\nu(x)$ is a state such that $\varphi_x(\lambda_g)=0$ for $\nu$-a.e.~$x\in X$ and every $g\notin G_x$. Consider the GNS-triple $(H,\xi,\pi)$ defined by $\Phi$. We claim that for $f,f_1,f_2\in C(X)$ and $g,h\in G$ we have
\begin{equation} \label{escalarpr}
(\pi(f(g\,\cdot)f_1\otimes\lambda_g)\xi, \pi(f_2\otimes\lambda_h)\xi)
=(\pi(f_1\otimes\lambda_g)\xi,\pi(\bar f(h\,\cdot)f_2\otimes\lambda_h)\xi).
\end{equation}
In other words, we claim that
$$
\int_Xf(gx)f_1(x)\bar f_2(x)\varphi_x(\lambda_{h^{-1}g})d\nu(x)=\int_Xf(hx)f_1(x)\bar f_2(x)\varphi_x(\lambda_{h^{-1}g})d\nu(x).
$$
Since $\varphi_x(\lambda_{h^{-1}g})=0$ for $\nu$-a.e.~$x$ not fixed by $h^{-1}g$, both integrals are in fact taken over the set $\{x\mid gx=hx\}$, but then they clearly coincide.

Similarly to the proof of Lemma~\ref{lcent}, it follows that we can define a representation $\rho$ of~$C(X)$ on~$H$ such that
$$
\rho(f)\pi(f_1\otimes\lambda_g)\xi=\pi(f(g\,\cdot)f_1\otimes\lambda_g)\xi\ \ \hbox{for}\ \ f,f_1\in C(X)\ \ \hbox{and}\ \ g\in G.
$$
Letting also $\rho(u_g)=\pi(1\otimes\lambda_g)$, we get a representation $\rho$ of $C(X)\rtimes G$ on $H$
such that $\rho(u_gf)\xi=\pi(f\otimes\lambda_g)\xi$. Then $\varphi:=(\rho(\cdot)\xi,\xi)$ is a state on $C(X)\rtimes G$ such that
$$
\varphi(u_gf)=\Phi(f\otimes\lambda_g)\ \ \hbox{for}\ \ f\in C(X)\ \ \hbox{and}\ \ g\in G.
$$
It is left to check that the centralizer of $\varphi$ contains $C(X)$, that is,
$$
\varphi(fu_gf_1)=\varphi(u_gff_1)\ \ \hbox{for}\ \ f,f_1\in C(X)\ \ \hbox{and}\ \ g\in G.
$$
As $fu_gf_1=u_gf(g\,\cdot)f_1$, this is a particular case of~\eqref{escalarpr} (with $f_2\equiv1$ and $h=e$).
\ep

The condition on $\Phi$ can be formulated as $\Phi(f\otimes\lambda_g)=0$ if $\supp f\cap X_g=\emptyset$, where $X_g=\{x\mid gx=x\}$. In this form the result is true without any separability assumptions on~$X$ and~$G$.

\smallskip

The condition $\varphi_x(\lambda_g)=0$ for $g\notin G_x$ implies that $\varphi_x$ is determined by its restriction to $C^*(G_x)$ and puts no conditions on this restriction. In other words, if $H$ is a subgroup of $G$ and $\psi\in S(C^*(H))$ is a state on $C^*(H)$, then there exists a unique state $\omega$ on $C^*(G)$ extending $\psi$ and such that $\omega(\lambda_g)=0$ for~$g\notin H$. Indeed, let $(H_\psi,\xi_\psi,\pi_\psi)$ be the GNS-triple defined by $\psi$. Consider the induced representation $\pi=\Ind^G_H\pi_\psi$ of $C^*(G)$. Recall that $\pi_\psi$ can be canonically identified with a subrepresentation of~$\pi|_{C^*(H)}$. Then $\omega=(\pi(\cdot)\xi,\xi)$. We will denote $\omega$ by $\Ind^G_H\psi$. Therefore Theorem~\ref{tmain} says that a state on $C(X)\rtimes G$ with centralizer containing $C(X)$ is determined by a probability measure $\nu$ on $X$ and a field of states $X\ni x\mapsto\psi_x\in S(C^*(G_x))$ defined $\nu$-a.e., with the only requirement that the map $(X,\nu)\to S(C^*(G))$, $x\mapsto \Ind^G_{G_x}\psi_x$, is weakly$^*$ measurable.

Note that induction also gives another way of constructing $\varphi$ out of $\Phi$. Namely, for every $x\in X$ the GNS-representation of $C^*(G_x)$ defined by $\varphi_x|_{C^*(G_x)}$ extends to a representation of $C(X)\rtimes G_x$ which maps $f\in C(X)$ into $f(x)1$. Hence there exists a state $\theta_x$ on $C(X)\rtimes G_x$ such that $\theta_x(fu_g)=f(x)\varphi_x(\lambda_g)$. Induce $\theta_x$ to a state~$\omega_x$ on~$C(X)\rtimes G$. Then $\varphi=\int_X\omega_xd\nu(x)$.

\smallskip

The subset of tracial states can be characterized as follows.

\begin{corollary}
Identity~\eqref{etrace} defines a bijection between tracial states $\varphi$ on $C(X)\rtimes G$ and states $\Phi=\int^\oplus_X\varphi_xd\nu(x)$ on $C(X)\otimes C^*(G)$ such that
\enu{i} $\Phi$ is invariant with respect to the action $\alpha$ of $G$ defined by $\alpha_g(f\otimes a)=f(g^{-1}\cdot)\otimes \lambda_ga\lambda_g^*$; equivalently, $\nu$ is $G$-invariant and $\varphi_{x}=\varphi_{gx}\circ\Ad\lambda_g$ for $\nu$-a.e. $x\in X$ and every $g\in G$;
\enu{ii} $\varphi_x(\lambda_g)=0$ for $\nu$-a.e. $x\in X$ and every $g\notin G_x$.
\end{corollary}

\bp If $\varphi$ is a state on $C(X)\rtimes G$ with centralizer containing $C(X)$, then $\varphi$ is tracial if and only if
$$
\varphi(u_gfu_hu_g^*)=\varphi(fu_h)\ \ \hbox{for}\ \ f\in C(X)\ \ \hbox{and}\ \ g,h\in G.
$$
Let $\Phi$ be the state on $C(X)\rtimes G$ corresponding to $\varphi$. Since $$\varphi(u_gfu_hu_g^*)=\varphi(f(g^{-1}\cdot)u_{ghg^{-1}})=\Phi(\alpha_g(f\otimes\lambda_h)),$$ it follows that $\varphi$ is a trace if and only if~$\Phi$ is $\alpha$-invariant.
\ep

Equivalently, we can say that to define a tracial state we need a $G$-invariant probability measure~$\nu$ on $X$ and a field of tracial states $\tau_x$ on $C^*(G_x)$ such that $\tau_{gx}=\tau_x\circ (\Ad\lambda_{g})^{-1}$ (so for every orbit $Gx\subset X$ the traces $\tau_y$ for $y\in Gx$ are determined by one trace $\tau_x$) and such that the map $(X,\nu)\ni x\mapsto\Ind^G_{G_x}\tau_x\in S(C^*(G))$ is weakly$^*$ measurable.

\smallskip

Now we move to the case when $G$ is abelian. For a subgroup $H$ of $G$ denote by $H^\perp$ the annihilator of $H$ in $\hat G$. For a measure $\lambda$ on $\hat G$ we have $\int_{\hat G}\chi(g)d\lambda(\chi)=0$ for all $g\notin H$ if and only if~$\lambda$ is $H^\perp$-invariant. Therefore we get the following result.

\begin{corollary} \label{cabel}
If $G$ is abelian, there is a bijection between tracial states $\tau$ on $C(X)\rtimes G$ and probability measures $\mu=\int^\oplus_X\nu_xd\nu(x)$ on~$X\times\hat G$ such that
\enu{i} $\mu$ is invariant with respect to the action of $G$ on the first factor of $X\times\hat G$; equivalently, $\nu$ is $G$-invariant and $\nu_{x}=\nu_{gx}$ for $\nu$-a.e. $x\in X$ and every $g\in G$;
\enu{ii} $\nu_x$ is $G_x^\perp$-invariant for $\nu$-a.e.~$x\in X$.

\noindent
Namely, the trace $\tau$ corresponding to such a measure $\mu$ is given by
$$
\tau(fu_g)=\int_{X\times\hat G}f(x)\chi(g)d\mu(x,\chi) \ \ \hbox{for}\ \ f\in C(X)\ \ \hbox{and}\ \ g\in G.
$$
\end{corollary}

Note that in this case by swapping the roles of $C(X)$ and $C(\hat G)$ we get yet another way of constructing the trace $\tau$ corresponding to $\mu$: define a representation~$\rho$ of $C(X)\rtimes G$ on~$L^2(X\times\hat G,d\mu)$ by
$$
(\rho(f)\zeta)(x,\chi)=f(x)\zeta(x,\chi),\ \ (\rho(u_g)\zeta)(x,\chi)=\chi(g)\zeta(g^{-1}x,\chi)
$$
and consider the function $\xi\equiv1$ in $L^2(X\times\hat G,d\mu)$, then $\tau=(\rho(\cdot)\xi,\xi)$.

\smallskip

Extremal tracial states correspond to extremal probability measures on $X\times\hat G$ with properties~(i) and~(ii). It is easy to see that such measures are of the form $\mu=\nu\times\lambda$, where $\nu$ is an ergodic $G$-invariant probability measure, the stabilizer $G_x$ does not depend on $x$ on a subset of full measure, denote this common stabilizer by $H$, and $\lambda$ is an extremal $H^\perp$-invariant measure, i.e., a point in $\hat G/H^\perp\cong\hat H$.

\begin{corollary}
Assume $G$ is abelian and $\tau$ is an extremal tracial state on $C(X)\rtimes G$. Then there exists a subgroup $H\subset G$, a character $\chi\in \hat H$ and an ergodic $G$-invariant probability measure $\nu$ on~$X$ with $G_x=H$ for $\nu$-a.e.~$x\in X$ such that
$$
\tau(fu_g)=\begin{cases}\chi(g)\int_Xf(x)d\nu(x), &\hbox{if}\ \ g\in H,\\
0,&\hbox{otherwise.}\end{cases}
$$
Conversely, any such triple $(H,\chi,\nu)$ defines an extremal tracial state.
\end{corollary}

We want to sketch a direct proof of this corollary. Assume $\tau$ is an extremal tracial state. Let~$\nu$ be the $G$-invariant probability measure on $X$ defined by $\tau|_{C(X)}$. The weak operator closure of $\pi_\tau(C(X))$ in $M:=\pi_\tau(C(X)\rtimes G)''$ is isomorphic to $L^\infty(X,\nu)$, and under this isomorphism the unitaries $\pi_\tau(u_g)$ implement the action of $G$ on $L^\infty(X,\nu)$. Since $M$ is a factor, it follows that~$\nu$ is ergodic. Next, for every $g\in G$ the set $X_g=\{x\mid gx=x\}$ is $G$-invariant, so either $X_g$ or its complement has measure zero. In other words, on a subset of full measure we have $G_x=H$ for a subgroup $H\subset G$. Since $M$ is a factor, for $h\in H$ the element $\pi_\tau(u_h)$ is scalar, so $\pi_\tau(u_h)=\chi(h)1$ for a character $\chi\in\hat H$. Thus, if $h\in H$ then $\tau(fu_h)=\chi(h)\tau(f)=\chi(h)\int_Xf\,d\nu$, and if $g\in G\setminus H$ then, using that $\nu(X_g)=0$, we get $\tau(fu_g)=0$ by the first part of the proof of Theorem~\ref{tmain}.

Conversely, assume we have such a triple $(H,\chi,\nu)$. Choose a section $s\colon G/H\to G$ of the quotient map $G\to G/H$, and define a $\T$-valued $2$-cocycle $\omega$ on $G/H$ by $\omega(g_1,g_2)=\chi(s(g_1)s(g_2)s(g_1g_2)^{-1})$. Since the action of $H$ on $(X,\nu)$ is essentially trivial, we have an action of $G/H$ on $L^\infty(X,\nu)$. Consider the twisted crossed product $M=L^\infty(X,\nu)\rtimes_{\omega}(G/H)$ in the von Neumann algebra sense. Since the action of $G/H$ is ergodic and essentially free, $M$ is a finite factor. Let $\tilde\tau$ be the unique normal trace on~$M$. Define a $*$-homomorphism $\rho\colon C(X)\rtimes G\to M$ by $\rho(f)=f\in L^\infty(X,\nu)$, $\rho(u_g)=\chi(gs(\bar g)^{-1})u_{\bar g}$, where $\bar g$ is the image of $g$ in $G/H$. Then $\tau:=\tilde\tau\circ\rho$ is the required trace, and as $\pi_\tau(C(X)\rtimes G)''\cong L^\infty(X,\nu)\rtimes_{\omega}(G/H)$ is a factor, $\tau$ is extremal.

\smallskip

Finally, consider the case $G=\Z$. Denote by $T\colon X\to X$ the homeomorphism corresponding to~$1\in\Z$. For $n\ge0$ denote by $X_n\subset X$ the subset of points of period $n$ (so $X_0$ is the set of aperiodic points). Then any measure $\mu$ on $X\times\T$ with properties (i) and (ii) from Corollary~\ref{cabel} decomposes into a sum of measures satisfying the same properties and concentrated on $X_n\times\T$ for some $n$.

If $\nu$ is concentrated on $X_0$ then $\mu=\nu\times\lambda$, where $\lambda$ is the Lebesgue measure, and the corresponding trace is $\nu_*\circ E$, where $\nu_*$ is the state on $C(X)$ defined by $\nu$ and $E\colon C(X)\rtimes G\to C(X)$ is the canonical conditional expectation.

If $\nu$ is concentrated on $X_n$, $n\ge1$, then $\mu$ is a $\Z/n\Z\times\Z/n\Z$-invariant measure on $X_n\times\T$, where the second factor $\Z/n\Z$ acts on $\T$ by rotations. Consider the simplest case when $\nu$ is concentrated on the orbit of a point~$x$ of period $n$. Then $\nu=n^{-1}\sum^{n-1}_{k=0}\delta_{T^kx}$ and $\mu=\nu\times\lambda$, where $\lambda$ is invariant under the rotation by~$2\pi/n$ degrees (we will say that $\lambda$ is $n$-rotation invariant). The corresponding trace can be written as follows.

The $*$-homomorphism $\rho\colon C(X)\to C(\Z/n\Z)$, $\rho(f)(k)=f(T^kx)$, extends to a $*$-homomorphism $$\rho\colon C(X)\rtimes\Z\to C(\Z/n\Z)\rtimes\Z.$$ Passing to the dual groups we can identify $C(\Z/n\Z)\rtimes\Z$ with $C(\T)\rtimes\Z/n\Z$.  By composing the canonical conditional expectation $C(\T)\rtimes\Z/n\Z\to C(\T)=C^*(\Z)$ with $\rho$ we then get a conditional expectation
$$
E_x\colon C(X)\rtimes\Z\to C^*(\Z),\ \ E_x(fu^m)=\frac{1}{n}\left(\sum^{n-1}_{k=0}f(T^kx)\right)u^m,
$$
where $u=u_1\in C(X)\rtimes\Z$ is the canonical unitary implementing $T$, $ufu^*=f(T^{-1}\cdot)$. The measure~$\lambda$ defines a state $\lambda_*$ on $C^*(\Z)$ by
$
\lambda_*(u^m)=\int_\T z^md\lambda(z).
$
Then $\lambda_*\circ E_x$ is the required tracial state on~$C(X)\rtimes\Z$.

It follows that in the case when there are only countably many periodic orbits Corollary~\ref{cabel} for~$G=\Z$ can be formulated as follows.

\begin{corollary} \label{ctrace}
Assume a homeomorphism $T$ of $X$ has at most countably many periodic orbits $O_i$, $i\in I$. For every $i$ choose $x_i\in O_i$. Then any tracial state $\tau$ on $C(X)\rtimes\Z$ has a unique decomposition
$$
\tau=\nu_*\circ E+\sum_{i\in I}{\lambda_i}_*\circ E_{x_i},
$$
where $\nu$ is a $T$-invariant measure on $X$ such that $\nu(O_i)=0$ for every $i$, $\lambda_i$ is an $|O_i|$-rotation invariant measure on $\T$, and $\nu(X)+\sum_i\lambda_i(\T)=1$. Conversely, any such collection of measures~$\nu$, $\lambda_i$ defines a tracial state.
\end{corollary}

\end{document}